\documentclass[11pt,reqno]{amsart}

\usepackage{graphics,wrapfig, graphicx, mathrsfs,xcolor}

\usepackage{amsfonts,amsthm,amsmath,amssymb,amscd,mathrsfs,graphicx,xcolor,subfig,tikz}
\usepackage{graphics}
\usepackage{indentfirst,centernot}
\usepackage{cite}
\usepackage{bm, enumerate,xcolor}
\usepackage[dvips]{epsfig}
\usepackage{latexsym}
\usepackage{float}
\usepackage{mathtools}
\usepackage{amssymb,amscd,graphicx,xcolor,subfig,tikz}
\usepackage{graphics}
\usepackage{indentfirst}
\usepackage{bm, enumerate,xcolor}
\usepackage[dvips]{epsfig}
\usepackage{latexsym}
\usepackage{float}
\usepackage{amsmath}
\usepackage{amssymb}
\newtheorem{lemma}{Lemma}[section]

\newtheorem{theorem}{Theorem}[section]
\newtheorem{corollary}{Corollary}[section]

\newtheorem{remark}{Remark}[section]

\numberwithin{equation}{section} \numberwithin{equation}{section}
\numberwithin{example}{section} \numberwithin{remark}{section}
\numberwithin{figure}{section} \numberwithin{algorithm}{section}
\mathtoolsset{showonlyrefs}

\def\ba{\begin{array}}
\def\ea{\end{array}}
\def\bma{\left(\begin{matrix}}
\def\ema{\end{matrix}\right)}
\def\be{\begin{equation}}
\def\ee{\end{equation}}

\def\vu{{\bf u}}

\def\dfrac{\displaystyle\frac}

\begin{document}
\title{Sequential Exponential Lower Bounds for the Advection-Diffusion Equation with Shear Flows}

\author[Y. Huang]{Yupei Huang}
\address{Yupei Huang:  Imperial College London, Department of Mathematics, Exhibition Road, South Kensington, London SW7 2AZ. }
\email{y.huang1@imperial.ac.uk}

\author[X. Xu]{Xiaoqian Xu}
\address{Xiaoqian Xu: Duke Kunshan University, Zu Chongzhi Center, No. 8 Duke Avenue, Kunshan, Jiangsu Province, 215316, P.R. China}
\email{xiaoqian.xu@dukekunshan.edu.cn}

\date{\today}

\subjclass[2020]{}%

\keywords{}%

\thanks{}
\begin{abstract}
In this paper, we prove that the $L^2$ norm of mean-free solutions to the advection--diffusion equation on $\mathbb{T}^2$ with shear drifts satisfies a sequential \emph{exponential lower bound} in time. This lower bound shows that diffusion can fundamentally suppress passive-scalar mixing.
\end{abstract}
\maketitle

\section{Introduction}
We study the advection--diffusion equation in the two dimensional torus $\mathbb{T}^2=[-\pi,\pi]^2$.
\[
\partial_t \rho + \nabla\!\cdot(\mathbf{u}\rho) = \mu \Delta \rho.
\]
This is a fundamental model in mathematical physics that describes the evolution of a scalar quantity $\rho$---such as temperature, pollutant concentration, or chemical density---transported by an incompressible fluid \cite{isenberg1973heat,chatwin1985mathematical,zlatev1984implementation,amundson1950mathematics}. As discussed in \cite{danckwerts1952definition, miles2018diffusion}, two key processes govern the dynamics: \emph{filamentation} and \emph{homogenization}. \emph{Filamentation} refers to the creation of fine spatial structures in the scalar field through stretching and folding by the flow; it is associated with the advective term $\nabla\!\cdot(\mathbf{u}\rho)$ and corresponds to a reduction in the \emph{scale of segregation}. \emph{Homogenization} describes the dissipation-driven relaxation of the scalar toward its spatial average; it is associated with the diffusion term $\mu\Delta \rho$ and corresponds to a reduction in the \emph{intensity of segregation}.
In this paper, we focus on the case where $\mathbf{u}$ is a divergence-free vector field, so the velocity field is incompressible. In addition, we assume that $\rho$ has zero spatial mean.

One of the key goals in studying this model is to understand the interplay between
\emph{filamentation} and \emph{homogenization}. It is well known that
\emph{filamentation} can accelerate \emph{homogenization}: stretching and folding
create finer structures, which in turn enhance the action of diffusion. An active
research direction, known as \emph{enhanced dissipation}, has identified large
classes of velocity fields that significantly accelerate the decay of the scalar.
See, for instance, 
\cite{constantin2008diffusion,kiselev2016suppression,feng2019dissipation,
zelati2020relation,elgindi2025optimal,masmoudi2022stability,masmoudi2020enhanced,
bedrossian2016enhanced,bedrossian2017enhanced,wei2021transition,wei2020linear,
wei2021diffusion}.

However, less is understood about the opposite direction: the effect of
\emph{homogenization} on \emph{filamentation}. It was conjectured in
\cite{miles2018diffusion} that, although filamentation accelerates
homogenization, the diffusion-driven homogenization process may in turn
\emph{suppress} filamentation. The delicate balance between these two mechanisms
is intimately connected to the Batchelor scale \cite{batchelor1959small}, which
characterizes the smallest dynamically relevant length scale in the scalar field.
The primary goal of this paper is to provide a mathematical justification, in
a general setting, that \emph{homogenization can suppress filamentation}.

To mathematically quantify this interplay between advection and diffusion, we turn to negative Sobolev norms.
\begin{itemize}
\item In the \textbf{inviscid case}, \emph{filamentation} is quantified by the decay of negative Sobolev norms, where $\|\rho\|_{\dot{H}^{-1}}/\|\rho\|_{L^2}$ (or, equivalently, $\|\rho\|_{\dot{H}^{-1}}$ in the case of incompressible flow) serves as a mathematical measure of the characteristic mixing scale \cite{mathew2003multiscale,lin2011optimal,thiffeault2012using}. For velocity fields with uniformly bounded Sobolev regularity (i.e., $\mathbf{u} \in L_t^\infty W_x^{1,p}$), this scale decays at most exponentially \cite{seis2013maximal, iyer2014lower}.
\item In the \textbf{viscous regime}, the mathematical analysis must track how variance is transferred to smaller scales and then dissipated. For the advection-diffusion equation with a bounded divergence-free drift $\mathbf{u}$ on the torus, $\partial_{t}\rho + \nabla\cdot(\mathbf{u}\rho) = \mu\Delta \rho$, it is known that the $L^2$ norm has the following bound.
\begin{equation}\label{eqn:paraL2}
C_1 \exp(-C_2 e^{C_3 t}) \leq \|\rho(t,\cdot)\|_{L^2} \leq C_4 \exp(-C_5 t),
\end{equation}
for some constants depending on $\vu,\mu$ and the initial data \cite{poon1996qnique,hongjie2017,miles2018diffusion}. For completeness, we provide the proof of this estimate in the Appendix. The double-exponential lower bound arises because the ratio $\|\nabla\rho\|_{L^2} / \|\rho\|_{L^2}$, which controls the cascade rate, can grow at most exponentially in time. This ratio has the same scaling as $\|\rho\|_{L^2} / \|\rho\|_{\dot{H}^{-1}}$, a quantity central to the analysis of inviscid mixing.
\end{itemize}

Recent work has established that these lower theoretical bounds are sharp.

\begin{itemize}
\item For \textbf{inviscid mixing}, the exponential decay in time of $\|\rho\|_{\dot{H}^{-1}}$ for velocity fields with uniformly bounded Sobolev regularity is known to be sharp, with concrete examples in \cite{yao2017mixing,alberti2019exponential,elgindi2019universal,elgindi2025optimal}.
\item In the \textbf{viscous regime}, the double-exponential lower bound in \eqref{eqn:paraL2} has been shown to be sharp by \cite{rowan2025superexponential}, who constructed bounded velocity fields in 2D for which the $L^2$ norm of the solution decays at this extreme rate. This represents the most efficient possible collaboration between advection and diffusion.
\end{itemize}

In contrast to collaborative interactions, our paper shows a distinct and counter-intuitive dynamic: although filamentation can enhance dissipation, \textbf{diffusion can suppress mixing}. This suppression of mixing was observed numerically \cite{miles2018diffusion}.

More precisely, for shear flows $U(t,y) \in L^{\infty}_t L^{2}_y$ on the two-dimensional torus $\mathbb{T}^2$, the solution of 
\begin{equation}\label{eqn:sheard}
\partial_{t}\rho + U(t,y)\,\partial_x \rho = \mu\Delta \rho, 
\qquad 
\rho(0,x,y)=\rho_0(x,y),
\end{equation}
satisfies an \textbf{exponential lower bound}:
there exists a sequence $t_n \rightarrow \infty$ such that 
\[
\|\rho(t_n)\|_{L^2} \ge C_1 e^{-C_2 t_n},
\]
where the constants depend on $\mu,U$ and $\rho_0$. As a direct consequence,
\[
\limsup\limits_{t\rightarrow\infty}\frac{\|\rho(t)\|_{\dot{H}^{-1}}}{\|\rho(t)\|_{L^2}} >0.
\]
Thus, for such flows, the normalized mixing scale cannot converge to zero; it returns along a sequence to scales bounded away from zero.This behavior differs from well-known inviscid mixing by shear flows (see, for example, \cite{bedrossian2016enhanced}): when $\mu=0$, there exist shear flows and initial data for which the mixing scale decays algebraically to zero. The results in this paper demonstrate a genuine \textbf{suppression of mixing by diffusion}. For illustration, we also provide an explicit example of shear flows and initial data in Section \ref{chap2}, for which the mixing scale $\|\rho(t)\|_{\dot{H}^{-1}}/\|\rho(t)\|_{L^2}$ can be estimated explicitly.

This result naturally prompts an investigation into the mechanism by which  \emph{homogenization} suppresses \emph{filamentation}. To this end, we analyze three illustrative scenarios.

First, in Section~\ref{chap2}, we examine a concrete shear flow. Beyond the classical algebraic mixing results for shear flows (see~\cite{bedrossian2017enhanced}), we construct a bounded shear profile together with initial data such that the solution to  
\[
\partial_t \rho + U(t,y)\,\partial_x \rho = 0
\]
satisfies the exponential mixing estimate  
\[
\|\rho(t)\|_{\dot{H}^{-1}} \lesssim e^{-\beta t}.
\]
In other words, this shear flow is an exponential mixer in the inviscid setting, achieving the optimal decay of the $\dot{H}^{-1}$ norm that is characteristic of the uniformly bounded Sobolev regularity regime. When combined with the exponential lower bound we establish for the advection-diffusion equation, this example shows that diffusion can suppress mixing even for flows that achieve this optimal inviscid mixing rate: the viscosity prevents the mixing scale from decaying to zero.

In Section~\ref{chap3}, we present a pulsed-diffusion model that illustrates a contrast between sequential and simultaneous advection and diffusion; this demonstrates that a simple pulsed protocol can behave differently from the continuous-time setting of Theorem~\ref{main}. A pulsed--diffusion system was previously considered in~\cite{feng2019dissipation}. By analyzing such a system on $\mathbb{T}^2$, where intervals of pure shear advection alternate with intervals of pure diffusion, we prove that the mixing scale of the solution can decay to zero asymptotically. As a result, the solution decays super-exponentially in time. Compared with our main theorem,  this example demonstrates a key physical insight: the sequential application of advection and diffusion, as in the pulsed--diffusion model, cannot reproduce the inhibitory effect that arises when both mechanisms act simultaneously.

Finally, in Section~\ref{chap4}, we present a formal 1D toy model, motivated by an artificial three-dimensional gluing construction, that illustrates why naive pointwise estimates can be misleading. In this carefully designed scenario, a preliminary analysis might suggest the possibility of super-exponential decay. However, a rigorous proof shows that the $L^2$ norm decays only exponentially in time. This demonstrates that producing genuine super-exponential decay is fundamentally difficult, even when one is not constrained by physical boundary conditions.

\subsection*{Acknowledgments.} This work is partially supported by National Key R\&D Program of China (2021YFA1001200). Y. H. acknowledges partial funding from EPSRC Horizon Europe Guarantee EP/X020886/1 and NSF
DMS-2043024. X. X. thanks the Mathematics Department at Duke University for their hospitality during a visit, when part of this work was conducted. He was partially supported by Kunshan Shuangchuang Talent Program (kssc202102066). The authors
thank Tarek M. Elgindi and Alexander A. Kiselev for discussions. 

\subsection*{AI disclosure}
The authors used AI tools for language editing, formatting, and technical proofreading, including suggestions for correcting errors and improving proof presentation. All scientific content and final decisions are the authors' own work.

\section{Sequential Exponential Lower Bound for the Advection–Diffusion Equation.}\label{mainsec}
We begin by recalling a well-known regularity theorem for weak solutions to the advection--diffusion equation (see, for example, \cite{bonicatto2024weak}).

\begin{theorem}\label{thm:exist}
Consider the initial value problem on the torus $\mathbb{T}^2=[-\pi,\pi]^2$:
\begin{equation}\label{eqn:ad}
    \partial_t \rho + \nabla \!\cdot (\mathbf{u}\rho) = \mu\Delta \rho.
\end{equation}
Suppose that $\mathbf{u}$ is divergence-free, $\mathbf{u}\in L^{\infty}([0,T];L^2(\mathbb{T}^2))$, and that $\rho(0,\cdot)\in L^2(\mathbb{T}^2)$. Then there exists a unique solution $\rho$ to \eqref{eqn:ad} in the parabolic energy class such that 
\[
\rho \in L^\infty(0,T;L^2(\mathbb{T}^2)) \cap L^2(0,T;H^1(\mathbb{T}^2)) \cap C([0,T];L^2(\mathbb{T}^2)).
\]
Moreover, for all $T\in (0,\infty)$, we have \[
\|\rho(T,\cdot)\|_{L^2}^2=\| \rho(0,\cdot)\|_{L^2}^2-2\mu\int_{0}^{T} \|\nabla \rho(t,\cdot)\|_{L^2}^2dt.
\]
\end{theorem}

In this section, we discuss the global behavior of the advection–diffusion equation with shear drifts in $L_{t}^{\infty}L_{y}^{2}$. By scaling, we assume that the diffusion coefficient $\mu=1$.
Inspired by \cite{lax1956stability}, we establish the sequential exponential lower bound for the solution to the advection-diffusion equation.

\begin{theorem}\label{main}
Consider the advection-diffusion equation on $\mathbb{T}^2$:
\begin{equation}\label{shear}
    \partial_t \rho + U(t,y)\,\partial_x \rho = \Delta \rho.
\end{equation}
Assume that $\rho(0,\cdot) \in L^2_{x,y}$ is mean-zero and non-zero and that 
\[
\|U\|_{L^{\infty}_t L^2_y} < \infty.
\]
Then there exist a constant $c>0$, depending on $U$ and $\rho(0,\cdot)$, such that 
\begin{equation} 
   \limsup \limits_{t\rightarrow\infty} e^{ct}\|\rho(t,\cdot)\|_{L^2}>0
\end{equation}
\end{theorem}

\begin{remark}
    Based on a standard energy estimate,  $\|\rho(t,\cdot)\|_{L^2_{x,y}}$ admits an exponential upper bound.
\end{remark}
\begin{remark}\label{rmk:couette}
It is well known that if we instead consider the advection--diffusion equation on $\mathbb{T}\times \mathbb{R}$ with the Couette flow drift $U(t,y)=y$, then for mean-zero-in-$x$ initial data the solution to \eqref{shear} exhibits super-exponential decay of its $L^2$ norm in time (see, for example, \cite{kelvin1887stability,bedrossian2016enhanced,masmoudi2022stability}).
 The key difference between the Couette flow in this setting and the flows covered by the theorem above is that the Couette flow is neither bounded nor integrable. 
\end{remark}
We now present the proof of Theorem \ref{main}.
\begin{proof}[Proof of Theorem \ref{main}.]
Suppose, to the contrary, that the conclusion fails. Then there exists a nonzero solution that decays super-exponentially, in the sense that for all $\lambda>0$, there exists $C_\lambda>0$ such that
\[
\|\rho(t)\|_{L^2_{x,y}} \le C_\lambda e^{-2\lambda t}.
\] In particular, $e^{\lambda t}\rho \in L_{t}^{2}L_{x,y}^2(\mathbb{R}^{+}\times \mathbb{T}^2)$.  Taking the Fourier expansion of $\rho:$ $\rho(t,x,y)=\sum\limits_{k\in\mathbb{Z}}\rho_{k}(t,y)e^{ikx}$, we are going to show that $\rho_{k}=0$, for all $k\in\mathbb{Z}$. 
Now assume that there is a $k$ such that $\rho_{k}\neq 0$. We have the following equation governing $\rho_{k}$:
    \begin{equation}
        \partial_{t}\rho_{k}+(k^2-\partial_{yy})\rho_{k}=-ikU(t,y)\rho_{k}, \text{for all $t\geq 0$.}
    \end{equation}
    Take $\Lambda_{m}=k^2+\frac{m^2+(m+1)^2}{2}$ with $m\in \mathbb{N}$, and define $h(t,y)$ as follows:
\begin{equation}
\begin{aligned}
   &h(t,y)=e^{\Lambda_{m}t}\rho_{k}(t,y), \quad\text{if $t\geq 0$},\\&
   h(t,y)=e^{\Lambda_{m}t}\rho_{k}(0,y)=e^{\Lambda_{m}t}\int_{\mathbb{T}}\rho(0,x,y)e^{ikx}dx,\quad \text{if $t\leq 0$.}
\end{aligned}
    \end{equation}
     By the Plancherel theorem, 
\begin{equation}
    \|\rho(t,\cdot)\|^2_{L^2_{x,y}}=2\pi\sum\limits_{k}\|\rho_k(t,\cdot)\|^2_{L^2_y},
\end{equation} 
we have \begin{equation}
    \int_{\mathbb{R}}|h(t,\cdot)|_{L_y}^2dt\leq \frac{1}{2\pi}[\int_{\mathbb{R^{+}}}e^{2\Lambda_{m}t}|\rho(t,\cdot)|_{L^2_{x,y}}^2dt+\frac{1}{2\Lambda_{m}}|\rho(0,\cdot)|_{L_{x,y}^2}^2]<+\infty.
\end{equation}
Similarly, we have $h\in L_{t}^1L_{y}^2$.
  Moreover,by the Cauchy-Schwarz inequality, we have \[
\|U(t)h(t)\|_{L^1_y}
\le
\|U(t)\|_{L^2_y}\|h(t)\|_{L^2_y}.
\]
Since $U \in L^\infty_t L^2_y$ and $h\in L^1_t L^2_y$, we have $Uh \in L^1_t L^1_y$. Also, $h \in L^\infty_t L^2_y$ (by continuity and the decay estimate), so $Uh \in L^\infty_t L^1_y$. Hence
\[
U(t,y) h(t,y) 1_{t>0}
\in
L^{1}_{t}L_{y}^{1}(\mathbb{R}\times \mathbb{T})
\cap
L^{\infty}_tL^1_y(\mathbb{R}\times \mathbb{T}).
\]


Moreover, since $h$ is continuous in $L^2$ in time, it is a distribution solution to the following equation.
    \begin{equation}
        \partial_{t}h+(k^2-\partial_{yy}-\Lambda_{m})h=-ik U(t,y)h 1_{t>0}+(k^2-\partial_{yy})\rho_{k}(0,\cdot) e^{\Lambda_{m}t}1_{t<0}:=F_1+F_2.
    \end{equation}    
Now we write the Fourier expansions of $h$, $\rho_k(0,\cdot)$, $F_1$, and $F_2$ in the $y$-variable and denote the corresponding Fourier coefficients by $h^{\,l}(t)$, $g_k^{\,l}$, $F_1^{\,l}(t)$, and $F_2^{\,l}(t)$, respectively. We obtain the following identity.

     \begin{equation}
        \partial_{t}h^{l}+(k^2+l^2-\Lambda_{m})h^{l}=F_1^{l}+F_2^{l}=F_{1}^{l}+1_{t\leq 0} (k^2+l^2)e^{\Lambda_{m}t} g_{k}^{l}.
    \end{equation}
    Now we define the Fourier transform operator in time \begin{equation}
         \hat{F}(\tau):=\int_{\mathbb{R}} F(t)e^{i\tau t}dt,
     \end{equation} 
using the standard distribution theory \cite{MR1996773}, and we have 
\begin{equation}\label{eqn:distribution Fourier}
\begin{aligned}
\hat{h}^l(\tau)=(-i\tau +k^2+l^2-\Lambda_{m})^{-1}[\Hat{F^l_1}(\tau)]+(-i\tau +k^2+l^2-\Lambda_{m})^{-1}[\frac{1}{\Lambda_{m}+i\tau}(k^2+l^2)g^l_{k}]
\end{aligned}
\end{equation}
for $l\in\mathbb{Z}$.

The key element of our proof lies in the analysis of equation \eqref{eqn:distribution Fourier}, and we will make use of the following lemma. \emph{From now on, throughout this section, $C$ denotes a positive constant independent of $m$, whose value may change from line to line.}
\begin{lemma}\label{lemma1}
{For any $k,l\in \mathbb{Z}$, $m\in \mathbb{Z}^+$} and $\frac{1}{2}<s<1$, let  $\Lambda_{m}=k^2+\frac{m^2+(m+1)^2}{2}$, {there exists a constant $C>0$ depending on $s$ and independent of $k,l,m$, such that} we have the following inequalities.

\begin{equation}\label{eqn:claim0}
\left(l^2-\frac{m^2+(m+1)^2}{2}\right)^2\geq C m^2;
\end{equation}

\begin{equation}\label{eqn:claim1}
         |(i\tau +k^2+l^2-\Lambda_{m})^{-1}|\leq \frac{C}{\sqrt{m^2+\tau^2}};
\end{equation}
\begin{equation}\label{eqn:claim2}
 |(i\tau +k^2+l^2-\Lambda_{m})^{-1}|\leq  \frac{C}{m^{1-s}|l|^{s}},\text{for $l \neq 0$};
\end{equation}
and
\begin{equation}\label{eqn:claim3}
 |(i\tau +k^2-\Lambda_{m})^{-1}|\leq \frac{C}{\sqrt{\tau^2+m^4}}.
\end{equation}
    \end{lemma}
\begin{proof}[Proof of the Lemma \ref{lemma1}.]
When $|l|\geq m+1$, we have \begin{equation}
    l^2-\frac{m^2+(m+1)^2}{2}\geq \frac{(m+1)^2-m^2}{2}\geq m>0.
\end{equation}
When $|l|\leq m$, we have \begin{equation}
    \frac{m^2+(m+1)^2}{2}-l^2\geq \frac{(m+1)^2-m^2}{2}\geq m>0.
    \end{equation}
    By combining the above two estimates, we have proved 
\eqref{eqn:claim0}.

 Equations \eqref{eqn:claim1}, \eqref{eqn:claim3} are direct consequences of equation \eqref{eqn:claim0}.

 Now we prove \eqref{eqn:claim2}.

 In the case where $|l|\geq 2m$, we have $l^2-\frac{m^2+(m+1)^2}{2}\geq \frac{l^2}{4}\geq m^2.$ Hence,  \begin{equation}
 \begin{aligned}
      &\quad |(i\tau+k^2+l^2-\Lambda_{m})^{-1}|=\frac{1}{\sqrt{\tau^2+[l^2-\frac{m^2+(m+1)^2}{2}]^2}}\leq \frac{1}{|l^2-\frac{m^2+(m+1)^2}{2}|}\\&\leq \frac{1}{(\frac{1}{4})^{s}}\frac{1}{m^{1-s}|l|^{2s}}\leq \frac{4^{s}}{m^{1-s}|l|^{s}}
 \end{aligned}
 \end{equation}
  In the case where $|l|\leq 2m$, we have $|l^2-\frac{m^2+(m+1)^2}{2}|\geq m>\frac{|l|}{4}$, hence, we have 
\begin{equation}
     |(i\tau+k^2+l^2-\Lambda_{m})^{-1}|=\frac{1}{\sqrt{\tau^2+[l^2-\frac{m^2+(m+1)^2}{2}]^2}}\leq \frac{1}{|l^2-\frac{m^2+(m+1)^2}{2}|}\leq  4^{s}\frac{1}{m^{1-s}|l|^{s}}.
 \end{equation}
Now, by combining the two estimates, we finish the proof of \eqref{eqn:claim2}. 
\end{proof}

Now, by \eqref{eqn:claim3} and \eqref{eqn:claim2} respectively in Lemma \ref{lemma1}, we have the following estimate:
\begin{equation}\label{eqn:eq0}
         |(i\tau +k^2-\Lambda_{m})^{-1} \Hat{F^{0}_1}(\tau)|\leq \frac{C}{\sqrt{m^4+
         \tau^2}}|\Hat{F^{0}_1}(\tau)| ;
\end{equation}
%
\begin{equation}\label{eqn:neq0}
         |(i\tau +k^2+l^2-\Lambda_{m})^{-1} \Hat{F^l_1}(\tau)|\leq \frac{C}{m^{1-s}|l|^{s}}|\Hat{F^l_1}(\tau)| \text{,\, for $l\neq 0$;} 
\end{equation}
Meanwhile, by \eqref{eqn:claim0}, we have  {
     \begin{equation}\label{eqn:estimate low to high0}
        \begin{aligned}
            &\quad \left|(i\tau +k^2+l^2-\Lambda_{m})^{-1}[\frac{1}{\Lambda_{m}+i\tau}(k^2+l^2)g^l_{k}]\right|\\&\leq \frac{k^2+l^2}{\sqrt{(k^2+\frac{m^2+(m+1)^2}{2})^2+\tau^2}}\frac{1}{\sqrt{(l^2-\frac{m^2+(m+1)^2}{2})^2+\tau^2}}|g_k^l|\\&
            =\frac{k^2}{\sqrt{(k^2+\frac{m^2+(m+1)^2}{2})^2+\tau^2}}\frac{1}{\sqrt{(l^2-\frac{m^2+(m+1)^2}{2})^2+\tau^2}}|g_k^l|\\&+
            \frac{l^2}{\sqrt{(k^2+\frac{m^2+(m+1)^2}{2})^2+\tau^2}}\frac{1}{\sqrt{(l^2-\frac{m^2+(m+1)^2}{2})^2+\tau^2}}|g_k^l|\\&
            =:J_1+J_2.
%
%
        \end{aligned}
    \end{equation}}
    For $J_1$, by \eqref{eqn:claim0}, we have \begin{equation}\label{eqn:estimate low to high1}
        |J_1|\leq \frac{1}{\sqrt{(l^2-\frac{m^2+(m+1)^2}{2})^2+\tau^2}}|g_k^l|\leq \frac{C}{\sqrt{m^2+\tau^2}}|g_k^l|.
    \end{equation}
    For $J_2$, the estimate is similar to the proof of \eqref{eqn:claim2}.
    In the  case where $|l|\leq 2m$, we have \begin{equation}
        |J_2|=\frac{l^2}{\sqrt{(k^2+\frac{m^2+(m+1)^2}{2})^2+\tau^2}}\frac{1}{\sqrt{(l^2-\frac{m^2+(m+1)^2}{2})^2+\tau^2}}|g_k^l|\leq \frac{16}{\sqrt{m^2+\tau^2}}|g_k^l|
    \end{equation}

    In the case where $|l|\geq 2m$, we have \begin{equation}
        |J_2|=\frac{1}{\sqrt{(k^2+\frac{m^2+(m+1)^2}{2})^2+\tau^2}}\frac{l^2}{\sqrt{(l^2-\frac{m^2+(m+1)^2}{2})^2+\tau^2}}|g_k^l|\leq \frac{16}{\sqrt{m^2+\tau^2}}|g_k^l|
    \end{equation}

    To sum up, there exists a constant $C>0$, independent of $m,l,k$ such that \begin{equation}\label{eqn:estimate low to high2}
        |J_2|\leq  \frac{C}{\sqrt{m^2+\tau^2}}|g_k^l|.
    \end{equation} 
    Now, combining \eqref{eqn:estimate low to high0},\eqref{eqn:estimate low to high1},\eqref{eqn:estimate low to high2}, we have for all $l \in \mathbb{Z}, $\begin{equation}\label{eqn:low to High}
        \left|(i\tau +k^2+l^2-\Lambda_{m})^{-1}[\frac{1}{\Lambda_{m}+i\tau}(k^2+l^2)g^l_{k}]\right|\leq \frac{C}{\sqrt{m^2+\tau^2}}|g_k^l|.
    \end{equation} 
    Moreover, for $l=0$, we have a sharper estimate:   
     \begin{equation}\label{eqn:0mode}
        \begin{aligned}
            &\quad \left|(i\tau +k^2-\Lambda_{m})^{-1}[\frac{1}{\Lambda_{m}+i\tau}k^2g^0_{k}]\right|\\&\leq \frac{k^2}{\sqrt{(\frac{m^2+(m+1)^2}{2})^2+\tau^2}}\frac{1}{\sqrt{(\frac{m^2+(m+1)^2}{2})^2+\tau^2}}|g_k^0|\\ &\leq C|g_k^0|\left(\frac{k^2}{\sqrt{m^4+\tau^2}}\right).
%
%
        \end{aligned}
    \end{equation}
   {
Now, based on  \eqref{eqn:eq0}, \eqref{eqn:neq0}, \eqref{eqn:low to High} and \eqref{eqn:0mode}, we have 
\begin{equation}
\begin{aligned}
        &\|\hat{h}(\tau)\|_{L_{y}^2}^2=\sum\limits_{l}|\hat{h}^l(\tau)|^2\leq C[\frac{|\hat{F}^0_1(\tau)|^2}{m^{4}+\tau^2}+|g_k^0|^2\left(\frac{k^4}{m^4+\tau^2}\right)\\&+\sum\limits_{l\neq 0}\frac{|\hat{F}^l_1(\tau)|^2}{l^{2s} m^{2-2s}}+|g_k^l|^2\frac{1}{m^2+\tau^2}]
\end{aligned}
    \end{equation}
    }
    Using Plancherel’s identity, we have 
\begin{equation}\label{eqn:positive time}
    \begin{aligned}
          &\quad \int_{\mathbb{R}^{+}} \|h(t,\cdot)\|_{L_{y}^{2}}^2dt \leq \int_{\mathbb{R}} \|h(t,\cdot)\|_{L_{y}^{2}}^2dt\\&\leq \int_{\mathbb{R}}C[\frac{|\hat{F}^0_1(\tau)|^2}{m^{4}+\tau^2}+|g_k^0|^2\left(\frac{k^4}{m^4+\tau^2}\right)+\sum\limits_{l\neq 0}\frac{|\hat{F}^l_1(\tau)|^2}{l^{2s} m^{2-2s}}+|g_k^l|^2\frac{1}{m^2+\tau^2}]d\tau\\&\leq \int_{\mathbb{R}}C[\frac{|\hat{F}^0_1(\tau)|^2}{m^{4}}+|g_k^0|^2\left(\frac{k^4}{m^4+\tau^2}\right)+\sum\limits_{l\neq 0}\frac{|\hat{F}^l_1(\tau)|^2}{l^{2s} m^{2-2s}}+|g_k^l|^2\frac{1}{m^2+\tau^2}]d\tau\\&\leq \frac{C\int_{\mathbb{R}}\|F_1(t,\cdot)\|_{H^{-s}_{y}}^2dt}{m^{2-2s}}+\frac{C}{m^4}\int_{\mathbb{R}}|F^0_{1}(t)|^2dt+C\int_{\mathbb{R}}\left(\frac{k^{4}}{m^4+\tau^2}+\frac{1}{m^2+\tau^2}\right)d\tau\|\rho_{k}(0,\cdot)\|_{L_{y}^2}^2.
    \end{aligned}
    \end{equation}
To close the argument, we note that for all $t>0$, $s>\frac{1}{2}$, by Sobolev embedding, we have:
    \begin{equation}\label{eqn:L^2 1}
        \|F_{1}(t,\cdot)\|_{H_{y}^{-s}}=\sup_{\|\phi\|_{H_{y}^{s}}=1} \langle F_{1}, \phi \rangle_{L_{y}^2}=\sup_{\|\phi\|_{H_{y}^{s}}=1} \langle ik U(t,y)h,\phi\rangle_{L_{y}^2} \leq |k|\|U(t,\cdot)\|_{L_{y}^2}\|h\|_{L_{y}^2}
    \end{equation}
    and \begin{equation}\label{eqn:L^2 2}
        |F^0_1(t)|=|\int_{\mathbb{T}} ik U(t,y)h(y)dy|\leq |k|\|U(t,\cdot)\|_{L_{y}^2}\|h\|_{L_{y}^2} 
    \end{equation}
    Combining \eqref{eqn:positive time}, \eqref{eqn:L^2 1} and \eqref{eqn:L^2 2}, we obtain
    \begin{equation}
    \begin{aligned}
         &\quad \int_{\mathbb{R}^{+}} \|h(t,\cdot)\|_{L_{y}^{2}}^2dt\leq \frac{Ck^2}{m^{2-2s}}\|U\|^2_{L_{t}^{\infty}L_{y}^2} \int_{\mathbb{R}^{+}} \|h(t,\cdot)\|_{L_{y}^{2}}^2dt+ \frac{Ck^2\|U\|^2_{L_{t}^{\infty}L_{y}^2}}{m^{4}}\int_{\mathbb{R}^{+}} \|h(t,\cdot)\|_{L_{y}^{2}}^2dt\\&+ C\|g_{k}\|_{L_{y}^2}^2 \int_{\mathbb{R}}\left(\frac{k^{4}}{m^4+\tau^2}+\frac{1}{m^2+\tau^2}\right)d\tau,
    \end{aligned}
    \end{equation}
    thus, for sufficiently large $m$, we have \begin{equation}
\begin{aligned}
      &\int_{\mathbb{R}^{+}}\|\rho_{k}(t,\cdot)\|_{L_{y}^{2}}^2dt \leq\int_{\mathbb{R}^{+}} e^{2(k^2+\frac{m^2+(m+1)^2}{2})t}\|\rho_{k}(t,\cdot)\|_{L_{y}^{2}}^2dt\\& =\int_{\mathbb{R}^{+}} \|h(t,\cdot)\|_{L_{y}^{2}}^2dt\leq \dfrac{C\|g_{k}\|_{L_{y}^2}^2 \left(\frac{k^4}{m^2}+\frac{1}{m}\right)}{1-\frac{Ck^2\|U\|^2_{L_{t}^{\infty}L_{y}^2}}{m^{2-2s}}-\frac{Ck^2\|U\|^2_{L_{t}^{\infty}L_{y}^2}}{m^{4}}}
\end{aligned}
    \end{equation}
Letting $m\rightarrow +\infty$, we have  $\int_{0}^{\infty}\|\rho_{k}(t,\cdot)\|_{L_{y}^{2}}^2dt = 0$, so $\rho_{k}(t,\cdot)=0$ for almost every $t>0$. Since $\rho_k \in C([0,\infty); L^2_y)$, we have $\rho_k(t,\cdot)=0$ for all $t\geq 0$. This contradicts the assumption $\rho_k \not\equiv 0$. Therefore $\rho_k=0$ for every $k\in\mathbb Z$, and hence $\rho\equiv 0$. 
\end{proof}
\begin{remark}\label{rmk:3d}
The same argument extends to the three-dimensional shear flow $(U(y,z,t),0,0)$ on $\mathbb{T}^3$, provided $U \in L^\infty_{t,y,z}$. We sketch the proof here.

In the Fourier variables $(k,\ell)$ with $k\in \mathbb Z$ and $\ell=(\ell_1,\ell_2)\in \mathbb Z^2$, the equation for the $k$-th Fourier mode is
\[
\partial_t \rho_k + (k^2 - \Delta_{y,z})\rho_k = -ik U(t,y,z)\rho_k.
\]
The resolvent estimate in Lemma 2.1 is replaced by a spectral gap argument in the two-dimensional symbol space. Let
\[
\mathcal S := \{ |\ell|^2 : \ell \in \mathbb Z^2 \} \subset \mathbb Z_{\ge 0}.
\]
By the Landau-Ramanujan theorem, the number of elements of $\mathcal S$ up to $N$ satisfies
\[
\#(\mathcal S \cap [0,N]) \sim K \frac{N}{\sqrt{\log N}}.
\]
The number $K$ is called Landau-Ramanujan constant. As a result, $\mathcal S$ has zero density and contains arbitrarily large consecutive gaps. Choose $\alpha_j \to \infty$ such that
\[
d_j := \mathrm{dist}(\alpha_j, \mathcal S) \to \infty.
\]
Set $\Lambda_j := k^2 + \alpha_j$. Then for every $\ell \in \mathbb Z^2$,
\[
|k^2 + |\ell|^2 - \Lambda_j| = \bigl||\ell|^2 - \alpha_j\bigr| \ge d_j,
\]
so the resolvent satisfies
\[
\left\| (-i\tau + k^2 + |\ell|^2 - \Lambda_j)^{-1} \right\|_{L^2 \to L^2}
\le \frac{1}{\sqrt{\tau^2 + d_j^2}}.
\]
The forcing term is bounded by
\[
\| -ik U h \|_{L^2(\mathbb T^3)} \le |k| \|U\|_{L^\infty_{t,y,z}} \|h\|_{L^2(\mathbb T^3)}.
\]
Thus the absorption coefficient becomes
\[
\frac{C k^2 \|U\|_{L^\infty}^2}{d_j^2},
\]
which tends to zero as $j\to\infty$ since $d_j\to\infty$. The initial-data term is bounded by an integral of the form
\[
\int_{\mathbb R} \left|\frac{k^2 + |\ell|^2}{(\Lambda_j + i\tau)(|\ell|^2 - \alpha_j - i\tau)}\right|^2 d\tau
\lesssim \frac{1}{d_j} + \frac{1}{\Lambda_j},
\]
which also tends to zero as $j\to\infty$. The same limiting argument as in the proof of Theorem 2.2 then yields the exponential sequential lower bound. We omit the routine details.
\end{remark}
Theorem 2.2 yields the following corollary:
\begin{corollary}\label{cor:mix}
With the same assumption as in Theorem \ref{main}, if the solution $\rho$ is smooth, we must have 
\begin{equation}
\limsup\limits_{t\rightarrow\infty}\frac{\|\rho(t)\|_{\dot{H}^{-1}}}{\|\rho(t)\|_{L^2}}>0.
\end{equation}
\end{corollary}
\begin{proof}
We argue by contradiction. Assume otherwise,  
\begin{equation}
\lim\limits_{t\rightarrow\infty}\frac{\|\rho(t)\|_{\dot{H}^{-1}}}{\|\rho(t)\|_{L^2}}=0.
\end{equation}
By Cauchy's inequality, we have
\begin{equation}
\frac{\|\nabla\rho(t)\|_{L^{2}}}{\|\rho(t)\|_{L^2}}\geq \frac{\|\rho(t)\|_{L^2}}{\|\rho(t)\|_{\dot{H}^{-1}}}\rightarrow \infty.
\end{equation}
This means for any $M>0$, there is a $T>0$, so that $\|\nabla\rho(t)\|_{L^{2}}\geq M\|\rho(t)\|_{L^2}$ for $t\geq T$. As a consequence, multiplying both sides of \eqref{shear} by $\rho$ and integration by parts, one gets
\begin{equation}
\partial_t\|\rho\|_{L^2}^2=-2\|\nabla\rho\|_{L^2}^2\leq -2M^2\|\rho\|_{L^2}^2.
\end{equation}
By Gr\"{o}nwall's inequality and the arbitrariness of $M$, this means the decay rate of $\|\rho\|_{L^2}$ can be arbitrarily large, which is a contradiction.
\end{proof}

\section{Exponential mixing by a shear flow}\label{chap2}
In this section,  we construct a shear flow $(U(t,y),0)$ with $\|U(t,y)\|_{L_{t,y}^{\infty
}}<1$ and find an initial datum such that the solution to the passive scalar equation associated with the shear flow and the initial data has the exponential decay for $\dot{H}^{-1}$ norm.
\begin{theorem}\label{thm:passive grow}
    For any smooth functions $f_0(y)$ and $g_0(y)$, not both identically zero, there exist constants $C_0>0,C_1>0$ and $\beta_1>0,\beta_2>0$ depending only on $f_0$ and $g_0$, and a $U(t,y)$ smooth in $y$ with $\|U\|_{L^{\infty}}\leq 1$ and $\|U(t,\cdot)\|_{H_{y}^1}\leq C_0e ^{\beta_1 \lceil t\rceil^2}$, such that the following is true.

    The solution $\rho$ of the following passive scalar equation \begin{equation}
    \begin{aligned}
                &\partial_{t} \rho+ U(t,y)\partial_{x}\rho=0,\\&
                \rho(0,x,y)=f_0(y)\sin{x}+g_0(y)\cos{x},
    \end{aligned}
    \end{equation}
    satisfies $\|\rho(t,\cdot)\|_{\dot{H}^{-1}}\leq C_1 e^{-\beta_2t}$. Here $\lceil x\rceil$ denotes the ceiling function.
\end{theorem}

Here, we briefly outline the main idea behind the construction of $U(t,y)$.  
For each interval $[n,n+1]$, and for a fixed explicit constant $A$ with $|A|<1$, we choose a sufficiently large $k(n)\in\mathbb{N}$ and show that there exists a time $t(n)\in(0,1)$ such that, if on the subinterval $t\in[n,n+t(n)]$ we take the velocity field to be $(\sin(k(n)y),0)$, then
\[
\|\rho(n+t(n),\cdot)\|_{\dot{H}^{-1}}
    \;\leq\;
    \frac{|A|+2}{3}\,\|\rho(n,\cdot)\|_{\dot{H}^{-1}}.
\]
In the remaining interval $[n+t(n),\,n+1)$, we simply choose $U(t,y)=0$.  
With this construction, it follows that
\[
\|\rho(n,\cdot)\|_{\dot{H}^{-1}}
    \;\leq\;
    \|\rho_0\|_{\dot{H}^{-1}}
    \left(\frac{|A|+2}{3}\right)^{n}.
\]

The reason for choosing such a velocity field is that, for all times $t\geq 0$, there exist functions $f(t,y)$ and $g(t,y)$ such that
\[
\rho(t,x,y)=f(t,y)\sin x + g(t,y)\cos x.
\]
We then use the following identity, which relates the $H^1$ norm of $\rho(t,\cdot)$ to the functions $f(t,\cdot)$ and $g(t,\cdot)$.

\begin{lemma}\label{lemma:H1 norm calculation}
For all $f, g \in H_y^1$, we have
\[
\|\,f(y)\sin x + g(y)\cos x\,\|_{H_{x,y}^1}^2
    = 2\pi \|f\|_{L^2}^2
    + 2\pi \|g\|_{L^2}^2
    + \pi \|f'\|_{L^2}^2
    + \pi \|g'\|_{L^2}^2.
\]
\end{lemma}

\begin{proof}
    We have \begin{equation}
        \begin{aligned}
            &\quad \|f(y)\sin{x}+g(y)\cos{x}\|_{H_{x,y}^1}^2\\&=\int_{\mathbb{T}^2} [f(y)\sin{x}+g(y)\cos{x}]^2+[f(y)\cos{x}-g(y)\sin{x}]^2+[f^{'}(y)\sin{x}+g^{'}(y)\cos{x}]^2dxdy.
        \end{aligned}
    \end{equation}
    Now, by Fubini's theorem, we have 
    \begin{equation}
        \begin{aligned}
            &\quad \|f(y)\sin{x}+g(y)\cos{x}\|_{H_{x,y}^1}^2\\&=\int_{\mathbb{T}^2} f^2(y)+g^2(y)+f'^{2}(y)\sin^2{x} +2 f^{'}(y)g^{'}(y)\sin{x}\cos{x}+g'^{2}(y)\cos^2{x}dxdy\\&=2\pi\|f\|_{L^2}^2+2\pi\|g\|_{L^2}^2+\pi\|f^{'}\|_{L^2} ^2+\pi\|g^{'}\|_{L^2}^2.
        \end{aligned}
    \end{equation}
\end{proof}
The velocity we constructed in Theorem \ref{thm:passive grow} will be a piecewise autonomous shear flow; hence we will need the following lemma for the time-independent shear flow.
\begin{lemma}
Consider the following initial value problem
 \begin{equation}
    \begin{aligned}
                &\partial_{t} \rho+ U(y)\partial_{x}\rho=0,\\&
                \rho(0,x,y)=f_0(y)\sin{x}+g_0(y)\cos{x}.
    \end{aligned}
    \end{equation}
    We have 
\begin{equation}
\begin{aligned}
    &\quad \|\rho(t,\cdot)\|_{H^1}^2= \|\rho(0,\cdot)\|_{H^1}^2-2 \pi t\int_{-\pi}^{\pi} U^{'}(y)[f_0(y)g_0^{'}(y)-f_0^{'}(y)g_0(y)]\\&+\pi t^2\int_{-\pi}^{\pi} (U^{'}(y)f_0(y))^2+(U^{'}(y)g_0(y))^2dy.
\end{aligned}
\end{equation}    
\end{lemma}
\begin{proof}
  We recall that
\[
\rho(t,x,y)=f_0(y)\sin\bigl(x-tU(y)\bigr)+g_0(y)\cos\bigl(x-tU(y)\bigr).
\]
Using trigonometric identities, we rewrite this expression in the form
\begin{equation}\label{eqn:transport-explicit}
\begin{aligned}
\rho(t,x,y)
&=\Bigl[f_0(y)\cos\!\bigl(tU(y)\bigr)+g_0(y)\sin\!\bigl(tU(y)\bigr)\Bigr]\sin x\\
&\quad+\Bigl[-f_0(y)\sin\!\bigl(tU(y)\bigr)+g_0(y)\cos\!\bigl(tU(y)\bigr)\Bigr]\cos x.
\end{aligned}
\end{equation}
For convenience, we define
\[
A_t(y):=f_0(y)\cos\!\bigl(tU(y)\bigr)+g_0(y)\sin\!\bigl(tU(y)\bigr),\qquad
B_t(y):=-f_0(y)\sin\!\bigl(tU(y)\bigr)+g_0(y)\cos\!\bigl(tU(y)\bigr),
\]
so that
\[
\rho(t,x,y)=A_t(y)\sin x+B_t(y)\cos x.
\]

By Lemma~\ref{lemma:H1 norm calculation}, the $H^1$-norm satisfies
\[
\|\rho(t,\cdot)\|_{H^1_{x,y}}^{2}
=2\pi\bigl(\|A_t\|_{L^2}^{2}+\|B_t\|_{L^2}^{2}\bigr)
+\pi\bigl(\|A_t'\|_{L^2}^{2}+\|B_t'\|_{L^2}^{2}\bigr).
\]
Differentiating $A_t$ and $B_t$, we obtain
\[
\begin{aligned}
A_t'(y)&=f_0'(y)\cos\!\bigl(tU(y)\bigr)+g_0'(y)\sin\!\bigl(tU(y)\bigr)
+tU'(y)\bigl[-f_0(y)\sin\!\bigl(tU(y)\bigr)+g_0(y)\cos\!\bigl(tU(y)\bigr)\bigr],\\
B_t'(y)&=-f_0'(y)\sin\!\bigl(tU(y)\bigr)+g_0'(y)\cos\!\bigl(tU(y)\bigr)
-tU'(y)\bigl[f_0(y)\cos\!\bigl(tU(y)\bigr)+g_0(y)\sin\!\bigl(tU(y)\bigr)\bigr].
\end{aligned}
\]

Direct expansion and cancellation of the mixed trigonometric terms yields
\begin{equation}\label{eq:H1-evolution}
\begin{aligned}
&\quad \|\rho(t,\cdot)\|_{H^1_{x,y}}^{2}
\\&=2\pi \|f_0\|_{L^2}^{2}+2\pi \|g_0\|_{L^2}^{2}
+\pi \|f_0'\|_{L^2}^{2}+\pi \|g_0'\|_{L^2}^{2}+\pi t^{2}\!\left(\|U' f_0\|_{L^2}^{2}+\|U' g_0\|_{L^2}^{2}\right)\\&-2\pi t\!\int_{-\pi}^{\pi}U'(y)\bigl[f_0(y)g_0'(y)-g_0(y)f_0'(y)\bigr]\,dy\\
&=\|\rho(0,\cdot)\|_{H^1}^{2}
+2\pi t\int_{-\pi}^{\pi} U'(y)\bigl[-f_0(y)g_0'(y)+f_0'(y)g_0(y)\bigr]\,dy \\
&+\pi t^2\int_{-\pi}^{\pi}\Bigl[(U'(y)f_0(y))^{2}+(U'(y)g_0(y))^{2}\Bigr]\,dy.
\end{aligned}
\end{equation}
\end{proof}
Based on the formula above and the Cauchy-Schwarz inequality, we have the following lemma.
\begin{lemma} \label{lemma: H1 growth upper bound}
    Let $k_0$ be a positive integer, $U(y)=\sin{{k_0 y}}$ and $\rho$ solve the following initial value problem
 \begin{equation}
    \begin{aligned}
                &\partial_{t} \rho+ U(y)\partial_{x}\rho=0,\\&
                \rho(0,x,y)=f_0(y)\sin{x}+g_0(y)\cos{x},
    \end{aligned}
    \end{equation}
    then $\|\rho(t,\cdot)\|_{H^1}^2\leq 2\|\rho(0,\cdot)\|_{H^1}^2+10\pi k_0^2 t^2 \|\rho(0,\cdot)\|_{L^2}^2$.
\end{lemma}
Now, we need another elementary lemma regarding the estimate of the stationary phase that appears in the formula of the solution.
\begin{lemma}\label{lem:elementary}
For all $k_0\in\mathbb{N}$ and $t\in(0,1]$, we have
\begin{equation}
\int_{-\pi}^{\pi}\sin\!\bigl(t\sin(k_0x)\bigr)\,dx = 0,
\end{equation}
and
\begin{equation}\label{eqn:zeromode}
\left|\frac{1}{2\pi}\int_{-\pi}^{\pi}\cos\!\bigl(t\sin(k_0x)\bigr)\,dx\right|:= a(t) < 1.
\end{equation}
\end{lemma}

\begin{proof}
Since $\sin\!\bigl(t\sin(k_0x)\bigr)$ is odd, we obtain
\[
\int_{-\pi}^{\pi}\sin\!\bigl(t\sin(k_0x)\bigr)\,dx = 0.
\]

For \eqref{eqn:zeromode}, use periodicity:
\begin{equation}\label{eqn: k0 independent}
\frac{1}{2\pi}\int_{-\pi}^{\pi}\cos\!\bigl(t\sin(k_0x)\bigr)\,dx
= \frac{1}{2k_0\pi}\int_{-k_0\pi}^{k_0\pi}\cos\!\bigl(t\sin x\bigr)\,dx
= \frac{1}{2\pi}\int_{-\pi}^{\pi}\cos\!\bigl(t\sin x\bigr)\,dx.
\end{equation}
Since $|\cos(\cdot)|\le 1$ with the equality $|\cos(t\sin x)|=1$ occurs only at isolated points on $[0,2\pi]$, the average is strictly less than $1$, completing the proof of \eqref{eqn:zeromode}.
\end{proof}

\bigskip

Moreover, we have the following lemma for the $\dot{H}^{-1}$ norm of $\rho$.

\begin{lemma}\label{lemma:H-1 decay}
Let
\[
A:=\frac{1}{2\pi}\int_{-\pi}^{\pi}\cos\!\bigl(\sin x\bigr)\,dx,
\qquad
\rho_0(x,y)=f_0(y)\sin x+g_0(y)\cos x\neq 0,
\]
\[
k_0=\left\lceil \frac{36}{1-|A|}\cdot\frac{\|\nabla\rho_0\|_{L^2}}{\|\rho_0\|_{\dot{H}^{-1}}}\right\rceil,
\qquad
U(y)=\sin(k_0y).
\]
Consider the initial value problem
\[
\partial_t\rho + U(y)\,\partial_x\rho = 0,
\qquad
\rho(0,x,y)=\rho_0(x,y).
\]

Then for all $t\in[0,1]$,
\begin{equation}\label{eqn:H-1 bound 1}
\|\rho(t,\cdot)\|_{\dot{H}^{-1}} \le 2\,\|\rho(0,\cdot)\|_{\dot{H}^{-1}},
\end{equation}
and
\begin{equation}\label{eqn:H-1 bound 2}
\|\rho(1,\cdot)\|_{\dot{H}^{-1}} \le \frac{2|A|+1}{3}\,\|\rho(0,\cdot)\|_{\dot{H}^{-1}}.
\end{equation}
Consequently, there exists $t_0\in(0,1)$ such that
\begin{equation}\label{eqn:ivt}
\|\rho(t_0,\cdot)\|_{\dot{H}^{-1}} = \frac{|A|+2}{3}\,\|\rho(0,\cdot)\|_{\dot{H}^{-1}}.
\end{equation}
\end{lemma}

\begin{remark}\label{rmk:A}
We have
\[
A=\frac{1}{2\pi}\int_{-\pi}^{\pi}\cos(\sin x)\,dx = J_0(1)\approx 0.7652,
\]
a positive constant strictly less than $1$, where $J_0$ is the Bessel function of the first kind.
\end{remark}

\begin{proof}[Proof of Lemma \ref{lemma:H-1 decay}]
We decompose
\begin{equation}
\begin{aligned}
\rho(0,\cdot)
&= \rho_0^{1} + \rho_0^{2} + \rho_0^{3} + \rho_0^{4} \\
&:= \sum_{|k|\le \lceil k_0/4\rceil} \widehat{f_0}(k)\, e^{iky}\sin x
   \;+\; \sum_{|k|\ge \lceil k_0/4\rceil+1} \widehat{f_0}(k)\, e^{iky}\sin x \\
&\quad+\; \sum_{|k|\le \lceil k_0/4\rceil} \widehat{g_0}(k)\, e^{iky}\cos x
   \;+\; \sum_{|k|\ge \lceil k_0/4\rceil+1} \widehat{g_0}(k)\, e^{iky}\cos x.
\end{aligned}
\end{equation}
Correspondingly, for $a = 1,2,3,4$, we define
\[
\rho^{a}(t,\cdot) := e^{-t\,U(y)\partial_x}\rho_0^{a}(\cdot).
\]
Since the transport operator $e^{-tU(y)\partial_x}$ is an $L^2$ isometry, for $a=2,4$,
\begin{equation}\label{eqn:high}
\|\rho^a(t,\cdot)\|_{\dot{H}^{-1}}
\le \|\rho^a(t,\cdot)\|_{L^2}
= \|\rho^a(0,\cdot)\|_{L^2}
\le \frac{\|\rho^a(0,\cdot)\|_{H^1}}{\lceil k_0/4\rceil}
\le \frac{1-|A|}{9}\|\rho(0,\cdot)\|_{\dot{H}^{-1}}.
\end{equation}
Now, for $\rho^1(t,\cdot)$ and $\rho^3(t,\cdot)$, let us denote $\sum_{|k|\leq \lceil \frac{k_0}{4}\rceil} \widehat{f_0}(k)e^{ik{y}}$ as $f_0^{low}(y)$ and $\sum_{|k|\leq \lceil \frac{k_0}{4}\rceil} \widehat{g_0}(k)e^{ik{y}}$ as $g^{low}_0(y)$, then from \eqref{eqn:transport-explicit} we have
\begin{equation}
\begin{aligned}
\rho^1+\rho^3&= f^{low}_0(y)\cos(t\sin(k_0y))\sin x-f_0^{low}(y) \sin(t\sin(k_0y))\cos x\\
&+g^{low}_0(y)\sin(t\sin(k_0y))\sin x+g^{low}_0(y)\cos(t\sin(k_0y))\cos x\\
&= f^{low}_1(t,y)\sin x+f^{low}_2(t,y)\cos x+g^{low}_1(t,y)\sin x+g^{low}_2(t,y)\cos x
\end{aligned}
\end{equation}
Notice that for any integer $k$, By symmetry and Lemma \ref{lem:elementary}, the following identities hold.
\begin{equation}
\begin{aligned}
&\frac{1}{2\pi}\int_0^{2\pi}e^{iky}\cos(t\sin(k_0y))dy=0,\quad \mbox{if $k$ is not divisible by $k_0;$} \\
&\frac{1}{2\pi}\int_0^{2\pi}\cos(t\sin(k_0y))dy=a(t);\\
&\frac{1}{2\pi}\int_0^{2\pi}e^{iky}\sin(t\sin(k_0y))dy=0, \quad \mbox{if $k$ is not divisible by $k_0.$} 
\end{aligned}
\end{equation}
Since $|\cos(\cdot)|\le 1$ and the equality $|\cos(t\sin x)|=1$ occurs only at isolated points, the average satisfies $|a(t)|<1$ for $0<t\le 1$.

Moreover, for $0\le t\le 1$, we have $a(t)=J_0(t)>0$, since the Bessel function $J_0$ is positive on $[0,1]$ (its first zero is at approximately $2.4048$). Thus $0<a(t)<1$ for $0<t\le 1$, and $a(0)=1$.

Taking the low-frequency projection $|k|\le \lceil \frac{k_0}{4}\rceil$, only the zero mode of $\cos(t\sin(k_0y))$ survives, whereas $\sin(t\sin(k_0y))$ has no zero mode. Thus the low-frequency part is $a(t)(\rho_0^1+\rho_0^3)$, and the remainder contains frequencies $|k|\ge \lceil \frac{k_0}{4}\rceil$. Hence
\begin{equation}\label{eqn:lowsplit}
\begin{aligned}
\rho^1+\rho^3&=\frac{1}{2\pi}\int_0^{2\pi}\cos(t\sin(k_0y))dy\left(\rho_0^1+\rho_0^3\right)\\
&+\sum\limits_{|k|\geq \lceil \frac{k_0}{4}\rceil}\left(\widehat{f^{low}_1}(t,k)\sin x+\widehat{f^{low}_2}(t,k)\cos x+\widehat{g^{low}_1}(t,k)\sin x+\widehat{g^{low}_2}(t,k)\cos x\right)e^{iky}\\
&=I+II.
\end{aligned}
\end{equation}
Since $\rho_0^1+\rho_0^3$ is obtained from $\rho_0$ by removing high frequencies, its $\dot{H}^{-1}$ norm cannot exceed that of $\rho_0$. Hence
\begin{equation}\label{eqn:low}
\|I\|_{\dot{H}^{-1}}\leq a(t)\|\rho(0,\cdot)\|_{\dot{H}^{-1}}. 
\end{equation}
Moreover,
\begin{equation}\label{eqn:high2}
\|II\|_{\dot{H}^{-1}}\leq \frac{1}{\lceil \frac{k_0}{4}\rceil}\|\rho^1+\rho^3\|_{L^2}=\frac{1}{\lceil \frac{k_0}{4}\rceil}\|\rho^1_0+\rho^3_0\|_{L^2}\leq \frac{1}{\lceil \frac{k_0}{4}\rceil}\|\rho(0,\cdot)\|_{H^1}\leq \frac{1-|A|}{9}\|\rho(0,\cdot)\|_{\dot{H}^{-1}}.
\end{equation}

Combining \eqref{eqn:high}, \eqref{eqn:low}, and \eqref{eqn:high2}, we have
\[
\|\rho(t,\cdot)\|_{\dot{H}^{-1}}
\le \left(\frac{1-|A|}{3}+a(t)\right)\|\rho(0,\cdot)\|_{\dot{H}^{-1}}.
\]
Since $a(1)=A$ and $\frac{1-|A|}{3}+a(t)<2$ for $t\in[0,1]$, this yields
\eqref{eqn:H-1 bound 1} and \eqref{eqn:H-1 bound 2}.  
Continuity of $\|\rho(t)\|_{\dot{H}^{-1}}$ in $t$ implies \eqref{eqn:ivt}.
\end{proof}
Now we are able to finish the proof of Theorem \ref{thm:passive grow}.

\begin{proof}[Proof of Theorem \ref{thm:passive grow}]
Consider the initial value problem
\begin{equation}\label{eqn:passive grow}
\begin{aligned}
&\partial_{t}\rho + U(t,y)\,\partial_{x}\rho = 0,\\&
\rho(0,x,y) = f_0(y)\sin x + g_0(y)\cos x.
\end{aligned}
\end{equation}

We construct the velocity field $U(t,y)$ iteratively on each interval $[n,n+1]$.  
Assume $U(t,y)$ has already been defined on $[0,n]$, so that $\rho(n,\cdot)$ is known.  
Set
\[
k(n):=\left\lceil\frac{36}{1-|A|}\frac{\|\nabla \rho(n,\cdot)\|_{L^2}}{\|\rho(n,\cdot)\|_{\dot{H}^{-1}}}\right\rceil,
\qquad
A=\frac{1}{2\pi}\!\int_{-\pi}^{\pi}\cos(\sin x)\,dx.
\]
By Lemma \ref{lemma:H-1 decay}, there exists $t(n)\in(0,1)$ such that
\[
\|e^{-t(n)\sin(k(n)y)\partial_x}\rho(n,\cdot)\|_{\dot{H}^{-1}}
=
\frac{2+|A|}{3}\,\|\rho(n,\cdot)\|_{\dot{H}^{-1}}.
\]
We define
\[
U(t,y)=\sin(k(n)y),\quad t\in[n,n+t(n)),
\qquad
U(t,y)=0,\quad t\in[n+t(n),n+1).
\]

Now the flow repeatedly injects high–frequency shear at carefully chosen time intervals, producing sustained phase mixing and a cascade of energy to finer scales, which drives the \(H^{1}\) growth while simultaneously enhancing decay in \(\dot{H}^{-1}\).

The corresponding solution to \eqref{eqn:passive grow} satisfies
\[
\|\rho(n+1,\cdot)\|_{\dot{H}^{-1}}
=
\frac{|A|+2}{3}\,\|\rho(n,\cdot)\|_{\dot{H}^{-1}}.
\]
Hence,
\[
\|\rho(n,\cdot)\|_{\dot{H}^{-1}}
=
\left(\frac{|A|+2}{3}\right)^n \|\rho(0,\cdot)\|_{\dot{H}^{-1}}.
\]

Moreover, for any $t\in[n,n+1]$, Lemma \ref{lemma:H-1 decay} yields
\begin{equation}\label{eqn: H-1 decaying}
\|\rho(t,\cdot)\|_{\dot{H}^{-1}}
\le 2\|\rho(n,\cdot)\|_{\dot{H}^{-1}}
\le 2\left(\frac{|A|+2}{3}\right)^{n}\|\rho(0,\cdot)\|_{\dot{H}^{-1}}
\le 2\left(\frac{|A|+2}{3}\right)^{t-1}\|\rho(0,\cdot)\|_{\dot{H}^{-1}}.
\end{equation}

Next, we estimate the growth of the $H^1$ norm.  
Let $\alpha_n:=\|\nabla \rho(n,\cdot)\|_{L^2}$.  
By Lemma \ref{lemma: H1 growth upper bound},
\[
\alpha_{n+1}^2
\le 2\alpha_n^2
+10\pi\,k(n)^2\,\|\rho(n,\cdot)\|_{L^2}^2
\le \left(2+C\Big(\tfrac{|A|+2}{3}\Big)^{-2n}\frac{\|\rho(0,\cdot)\|_{L^2}^2}{\|\rho(0,\cdot)\|_{\dot{H}^{-1}}^2}\right)\alpha_n^2.
\]
Thus, there exist constants $C(f_0,g_0)$ and $\beta>0$ such that
\[
\alpha_{n+1} \le C(f_0,g_0)e^{\beta n}\alpha_n,
\]
and hence, by induction,
\[
\alpha_n \le C(f_0,g_0)e^{3\beta n^2}\|\nabla\rho(0,\cdot)\|_{L^2}.
\]

For $t\in[0,1]$, using Lemma \ref{lemma: H1 growth upper bound} again as in \eqref{eqn: H-1 decaying}, we obtain
\[
\|\nabla\rho(n+t,\cdot)\|_{L^2}
\le C(f_0,g_0)e^{\beta_3(n+1)^2}\|\nabla \rho(0,\cdot)\|_{L^2}.
\]

Finally, by our construction,
\[
\|U(t,\cdot)\|_{H_y^1}
\le k(\lfloor t\rfloor)
=\frac{36}{1-|A|}\frac{\|\nabla \rho(\lfloor t\rfloor,\cdot)\|_{L^2}}{\|\rho(\lfloor t\rfloor,\cdot)\|_{\dot{H}^{-1}}}
\le C(f_0,g_0)e^{\beta_3(\lfloor t\rfloor+1)^2}.
\]
\end{proof}

\section{An illustrative pulsed-diffusion example}\label{chap3}
The following elementary example illustrates that a pulsed-diffusion iteration (similar to that considered in \cite{feng2019dissipation}) can produce super-exponential decay. This highlights a key distinction from the simultaneous advection-diffusion setting of Theorem~\ref{main}, where such decay is impossible.

To see it, consider the following concrete iterative model:
\begin{equation}
\label{pulseddiffusion}
\rho(n,X)=e^{\Delta}\rho(n-1,TX), \qquad X=(x,y)\in \mathbb{T}^2,
\end{equation}
where $e^{\Delta}$ is is the time-one solution operator for the heat equation and $T=\begin{bmatrix}
1 & -1\\
0 & 1
\end{bmatrix},$ with
\begin{equation}\label{initialdata}    
\rho(0,{x},y)=e^{i(x+y)}+e^{-i(x+y)}.
\end{equation}

One can see that \eqref{pulseddiffusion} is a pulsed-diffusion version of the advection-diffusion equation with the advection velocity given by the Couette flow $v=(y,0)^{T}$. By direct calculations and induction, one can get the following:
\begin{theorem}
For $\rho$ defined by \eqref{pulseddiffusion} and \eqref{initialdata}, the mixing scale can be estimated as:
\begin{align*}
\begin{split}
\frac{\|\rho(n,\cdot)\|_{\dot{H}^{-1}}}{\|\rho(n,\cdot)\|_{L^{2}}}&=\frac{1}{\sqrt{1+(n+1)^2}}\rightarrow 0 \mbox{  as  } n\rightarrow \infty.
\end{split}
\end{align*}
In addition,
\begin{align*}
\|\rho(n,\cdot)\|_{L^2}=2\sqrt{2}\pi\exp(\frac{-2n^3-9n^2-19n}{6}).
\end{align*}
\end{theorem}
\begin{proof}

 We have the explicit formula for $\rho(n,\cdot)$:
   \begin{align*}
\begin{split}
\rho(n,x,y)&=\exp{(-\sum_{k=1}^{n}(1+(k+1)^2))}\cdot \left(e^{i(x-(n+1)y)}+e^{-i(x-(n+1)y)}\right)\\
&=\exp(\frac{-2n^3-9n^2-19n}{6})\left(e^{i(x-(n+1)y)}+e^{-i(x-(n+1)y)}\right).
\end{split}
\end{align*} 
As a result, we have 
    \begin{align*}
\begin{split}
&\frac{\|\rho(n,\cdot)\|_{\dot{H}^{-1}}}{\|\rho(n,\cdot)\|_{L^{2}}}=\frac{1}{\sqrt{1+(n+1)^2}},\\
&\|\rho(n,\cdot)\|_{L^2}=2\sqrt{2}\pi\exp(\frac{-2n^3-9n^2-19n}{6}).
\end{split}
\end{align*}
\end{proof}
\begin{remark}
The decay of the mixing scale of $\rho(n,\cdot)$ to $0$  as $n \rightarrow \infty$ shows a key difference from the result in Corollary \ref{cor:mix}. Moreover, the $L^2$ norm of $\rho$ decays faster than exponentially in terms of $n$. 
\end{remark}
\begin{remark}
One may also notice that, the iteration \eqref{pulseddiffusion} corresponds, up to a sign, to advection by a sawtooth shear on the torus, rather than a smooth periodic Couette field. We present it only as an illustrative discrete model.
\end{remark}

\section{A 1D model equation}\label{chap4}
We now present a formal heuristic construction, intended only to motivate why naive guesses about super-exponential decay can be misleading. We study a 1D model for the potential faster decay of $L^2$ norm of the advection-diffusion equation, with an artificial boundary condition. Nevertheless, we can show that this model does not yield faster-than-exponential decay, even under this highly artificial boundary condition.
Let $\mathbb{T}^2$ be the 2D torus, and let $A=\begin{bmatrix}2 & 1\\ 1 & 1\end{bmatrix}$ be the cat map. Given $v=(0,0,-1)^{T}$, for $\rho$ defined on $\mathbb{T}^2\times \mathbb{R}$ satisfying the following advection-diffusion equation with boundary value condition:
\begin{equation}\label{1Dmodel}
\begin{split}
&\partial_{t}\rho+v\cdot\nabla \rho-\Delta \rho=0, \\ &\quad\rho(x,y,z+1,t)=\rho(A^{T}(x,y)^{T},z,t), \quad\int_{\mathbb{T}^2}\rho(x,y,z,t)dxdy=0.
\end{split}
\end{equation}
This gluing is only formal since $A$ is not an isometry of $\mathbb{T}^2$. One would expect the solution $\rho$ to decay faster than exponentially in time for the following reason:

Taking the Fourier transform in $x$ and $y$, we get
\begin{equation}
\partial_{t}\hat{\rho}-\partial_z\hat{\rho}=\partial^2_z\hat{\rho}-|k|^2\hat{\rho}
\end{equation}
with boundary condition $$\hat{\rho}(k,1)=\hat{\rho}(A^{-1}k,0).$$

Hence, if we fix $k_0\in\mathbb{Z}^2$, and define $\theta(z)=\hat{\rho}(k_0,z)$, by boundary condition we have \eqref{1Dmodel} becomes
\begin{equation}
\partial_{t}\theta -\partial_z \theta -\partial_z^2 \theta =-|A^{-\lfloor z \rfloor}k_0|^2 \theta .
\end{equation}
Here $\lfloor z \rfloor$ is the floor function.

As a consequence, one could expect the following toy model will share the same decay property as of $\theta$:
\begin{equation}\label{toymodel}
\partial_{t} l-\partial_z l-\partial_z^2 l=-e^{|z|}l.
\end{equation}
Substituting $y=z+t$, one gets
\begin{equation}
\partial_{t}l-\partial_{yy}l=-e^{|y-t|}l.
\end{equation}
Such an equation is a 1D heat equation with extra force, defined on the whole line. One may want to show its decaying faster than exponential, given the following classical result for the expression of the solution, which is similar to the Feynman–Kac formula (say, one can find it from Theorem $2.2$ of Chapter 1 in the book by Freidlin \cite{freidlin1985functional}:

$$
l(t,z)=\mathbb E_z\left[l_0(X_t)\exp\left(-\int_0^t e^{|X_s|}\,ds\right)\right],
\qquad X_s = z + s + \sqrt{2}B_s.
$$ 
Here $W_t^z$ is the standard Brownian motion. Meanwhile, by the formula on page $36$ of \cite{freidlin1985functional}, a formal substitution using the law of the iterated logarithm, one has
$$
W_t^z\approx \sqrt{2t\ln\ln t}\quad a.s,
$$
which would seem to suggest 
$$
\exp \left(-\int_0^t e^{|t-s+B_{t-s}^{z}|}\,ds\right)\leq \exp \left(\int_{\frac{t}{2}}^t e^{\frac{1}{2}(t-s)}\,ds\right)\leq\exp \left(-e^{\frac{t}{4}}\right)\quad a.s. 
$$
This would suggest that the solution to \eqref{1Dmodel} decay faster than exponential, if the preceding estimate were valid. 

However, this is only a heuristic: the law of the iterated logarithm is a subsequential statement and provides no uniform bound over the whole integral $\int_0^t e^{|W_s|}\,ds$. The expectation is controlled by paths with small values of this integral, not by the extreme fluctuations described by the law of the iterated logarithm. The precise spectral analysis below shows that the actual decay is only exponential.

\begin{theorem}\label{thm:1d toy case lower bound}
    Let $l$ be a solution to \eqref{toymodel} with $l(0,z)=l_{0}(z)\not\equiv 0$. If $l_{0}(z)\in L^2(\mathbb{R})$ and $l_{0}(z)e^{\frac{z}{2}}\in L^2(\mathbb{R})$, then there exists $C(l_{0})>0$ such that \begin{equation}
        \|l(t,\cdot)\|_{L^2(\mathbb{R})}\geq \frac{1}{C(l_0)}e^{-C(l_{0})t}.
    \end{equation}
\end{theorem}
\begin{remark}
Compared with the unbounded Couette flow on $\mathbb{T}\times\mathbb{R}$ discussed in Remark \ref{rmk:couette}, the flow we have here, $v=(0,0,-1)^{T}$ on $\mathbb{T}^2\times\mathbb{R}$, is bounded but not integrable.
\end{remark}
Motivated by classical Sturm–Liouville theory, we take the following choice of change of variable $l(t,z)=m(t,z)e^{Q(z)}$.
The equation \eqref{toymodel} is then transformed to \begin{equation}
    \begin{aligned}
         \partial_{t}m = \partial_{zz}m + (1+2Q^{'}(z))\partial_{z}m + \bigl(Q^{'}(z)+Q^{''}(z)+(Q^{'})^2(z)-e^{|z|}\bigr)m.
    \end{aligned}
\end{equation}
Setting $Q(z)=-\frac{z}{2}$ and we end with \begin{equation}\label{eqn:evolution self-adjoint}
    \partial_{t}m=\mathcal{L}m:=\partial_{zz}m -(\frac{1}{4}+e^{|z|})m.
\end{equation}
Since $\frac{1}{4}+e^{|z|}>\frac{1}{4}$ and $\lim\limits_{|z|\rightarrow \infty}\frac{1}{4}+e^{|z|}=+\infty$, the operator $\mathcal{L}$ enjoys the following property. 
\begin{lemma}
Let
\[
\mathscr{L} := -\partial_{zz} + \Big(\tfrac14 + e^{|z|}\Big)
\qquad\text{acting on } L^2(\mathbb{R}),
\]
with domain
\[
\mathrm{Dom}(\mathscr{L})
:= \Big\{\, m\in L^2(\mathbb{R}) : m,m'\in AC_{\mathrm{loc}}(\mathbb{R}),\;
-m'' + (\tfrac14 + e^{|z|})m \in L^2(\mathbb{R}) \,\Big\},
\]
here $AC_{\mathrm{loc}}(\mathbb{R})$ is the class of locally absolutely continuous functions on $\mathbb{R}$.

Then:
\begin{itemize}

\item[(i)]
The operator $\mathscr{L}$ is self-adjoint on $L^2(\mathbb{R})$ and has compact resolvent( see, for example, \cite{berezin2012schrodinger}.) 
In particular, its spectrum is purely discrete and consists of a strictly increasing sequence of real eigenvalues
\[
0 < \lambda_0 < \lambda_1 < \lambda_2 < \cdots,
\qquad 
\lambda_n \to +\infty \text{ as } n\to\infty.
\]

\item[(ii)]
Each eigenvalue $\lambda_n$ is simple (geometric and algebraic multiplicity equal to $1$); see \cite{berezin2012schrodinger}.

\item[(iii)]
Define the even/odd subspaces
\[
L^2_{\mathrm{even}}(\mathbb{R}) := \{f\in L^2(\mathbb{R}) : f(z)=f(-z)\},\qquad
L^2_{\mathrm{odd}}(\mathbb{R}) := \{f\in L^2(\mathbb{R}) : f(z)=-f(-z)\}.
\]
Then
\[
L^2(\mathbb{R}) = L^2_{\mathrm{even}}(\mathbb{R}) \oplus L^2_{\mathrm{odd}}(\mathbb{R}),
\]
and each subspace is invariant under $\mathscr{L}$. 
\end{itemize}
\end{lemma}

Now we are ready to give the proof of Theorem \ref{thm:1d toy case lower bound}.
\begin{proof}
By the spectral properties established above, let $\{\phi_n\}_{n=0}^\infty$ be an orthonormal eigenbasis of $\mathscr L$ with corresponding eigenvalues $0<\lambda_0<\lambda_1<\cdots$. Since $m_0=e^{z/2}l_0\in L^2(\mathbb R)$ and $m_0\not\equiv 0$, we have the expansion
\[
m(t,z)=\sum_{n=0}^\infty c_n e^{-\lambda_n t}\phi_n(z),
\qquad c_n=\langle m_0,\phi_n\rangle_{L^2}.
\]
Choose $n_*$ be the smallest index such that  $c_{n_*}\neq 0$, we have 
\[
m(t) = c_{n_*} e^{-\lambda_{n_*}t}\phi_{n_*} + o_{L^2}(e^{-\lambda_{n_*}t}).
\]
The potential $V(z)=\frac14+e^{|z|}$ is even, and the eigenvalues are simple, so each $\phi_n$ is either even or odd. Hence
\[
\int_{\mathbb R_-} |\phi_{n_*}(z)|^2\,dz = \frac12.
\]
Therefore, for all sufficiently large $t$,
\[
\int_{\mathbb R_-} |m(t,z)|^2\,dz
\ge \frac{|c_{n_*}|^2}{4} e^{-2\lambda_{n_*}t}.
\]
Since $l=e^{-z/2}m$ and $e^{-z}\ge 1$ on $\mathbb R_-$,
\[
\|l(t)\|_{L^2(\mathbb R)}^2
\ge \int_{\mathbb R_-} |l(t,z)|^2\,dz
= \int_{\mathbb R_-} e^{-z}|m(t,z)|^2\,dz
\ge \int_{\mathbb R_-} |m(t,z)|^2\,dz
\ge C e^{-2\lambda_{n_*}t}.
\]
Taking square roots gives, for all sufficiently large $t$,
\[
\|l(t)\|_{L^2(\mathbb R)} \ge C e^{-\lambda_{n_*}t}.
\]
For small \(t\), the bound follows from continuity of $t\mapsto \|l(t)\|_{L^2}$ and the fact that $l_0\not\equiv 0$. Adjusting the constant if necessary, we obtain the desired exponential lower bound for all $t\ge 0$.
\end{proof}

\section{Appendix}

In this section, we prove the following theorem, following ideas similar to those in \cite{miles2018diffusion,hongjie2017}.

\begin{theorem}\label{thm:para}
Let $\mathbf{u}(t,\mathbf{x}) \in L^{\infty}_{t,\mathbf{x}}$ be a two-dimensional divergence free vector field in $\mathbb{T}^2$. Suppose $\rho$ is the unique non-trivial solution to
\begin{equation}\label{eqn:para}
\partial_{t}\rho + \mathbf{u}\cdot\nabla\rho = \mu\Delta\rho,
\end{equation}
with $\rho(0,\cdot)\in L^{2}(\mathbb{T}^2)$ and
\[
\int_{\mathbb{T}^2}\rho(0,\mathbf{x})\,d\mathbf{x}=0.
\]
Then there exists a constant $C_1,C_2>0$ where $C_1$ depends  on $\mu$, $\|\mathbf{u}\|_{L^{\infty}}$, and the initial datum and $C_2$ depends on $\mu$, $\mathbf{u}$, and the initial datum such that
\begin{equation}\label{eqn:double}
\|\rho(t,\cdot)\|_{L^2}\ge C_2\exp\!\bigl(-2\mu\exp(C_1t)\bigr),
\qquad t\ge0.
\end{equation}
\end{theorem}

\begin{proof}[Proof of Theorem \ref{thm:para}]
By the energy identity in Theorem \ref{thm:exist}, we have
\[
2\mu\!\int_0^\infty \|\nabla\rho(t,\cdot)\|_{L^2}^2\,dt
\le \|\rho(0,\cdot)\|_{L^2}^2,
\qquad
2\mu\!\int_0^T \|\nabla\rho(t,\cdot)\|_{L^2}^2\,dt
= -\|\rho(T,\cdot)\|_{L^2}^2+\|\rho(0,\cdot)\|_{L^2}^2,
\]
for all $T\ge0$. In particular, there exists $T_0<\infty$ such that
\[
\rho(T_0,\cdot)\in H^1(\mathbb{T}^2), \qquad \rho(T_0,\cdot)\neq0.
\]
The fact that $\rho(T_0)\neq 0$ follows from the backward uniqueness property of the parabolic operator $\partial_t + \mathbf u\cdot\nabla - \mu\Delta$: if $\rho(T_0)=0$, then by backward uniqueness $\rho(0)=0$, contradicting the assumption. See, e.g., \cite{poon1996qnique}.

By replacing $\rho(0,\cdot)$ with $\rho(T_0,\cdot)$, we may assume from now on that  
\[
\rho(0,\cdot)\in H^1(\mathbb{T}^2), \qquad \rho(0,\cdot)\neq 0, \qquad \int_{\mathbb{T}^2}\rho(0,\mathbf{x})dx=0. 
\]

To justify the differentiation of $\|\nabla\rho\|^2$ in the following lemma, we note that the computation can be made rigorous by a standard approximation argument: one regularizes the velocity field $\mathbf u$, derives the estimate for smooth solutions, and passes to the limit. Equivalently, one may use the standard $L^2$-maximal regularity result for the heat equation, which gives $\rho_t,\Delta\rho\in L^2_{\mathrm{loc}}((0,\infty);L^2)$ for $H^1$ initial data. Both approaches yield the same differential inequality for the energy solution.

We first prove the following key differential inequality.

\begin{lemma}\label{lem:ratio}
Under the assumptions of Theorem~\ref{thm:para}, suppose additionally that $0\not\equiv \rho(0,\cdot)\in H^1(\mathbb{T}^2)$. Then there exists $C_1>0$ such that
\begin{equation}\label{eqn:ratio}
\frac{\|\nabla\rho(t,\cdot)\|_{L^2}}{\|\rho(t,\cdot)\|_{L^2}}
\;\le\;
C_1 e^{C_1t},
\qquad t\ge0.
\end{equation}
\end{lemma}

\begin{proof}[Proof of Lemma \ref{lem:ratio}]
Define
\[
H(t)=\frac{\|\nabla\rho(t,\cdot)\|_{L^2}^2}{\|\rho(t,\cdot)\|_{L^2}^2}.
\]
A direct computation gives
\begin{equation}\label{eqn:HoL}
\partial_tH(t)
=
\frac{1}{\|\rho\|_{L^2}^4}\left(
\partial_t\|\nabla\rho\|_{L^2}^2\,\|\rho\|_{L^2}^2
-
\partial_t\|\rho\|_{L^2}^2\,\|\nabla\rho\|_{L^2}^2
\right).
\end{equation}

We compute each term. Using \eqref{eqn:para} and integration by parts:
\begin{equation}\label{eqn:H1}
\|\nabla\rho\|_{L^2}^2
= \int\rho(-\Delta\rho)
= -\frac{1}{2\mu}\!\int\rho\,(\mathbf{u}\cdot\nabla\rho)\,dx
-
\frac{1}{2\mu}\!\int\rho(\partial_t\rho+\mu\Delta\rho)\,dx.
\end{equation}
Similarly,
\begin{equation}\label{eqn:L2}
\partial_t\|\rho\|_{L^2}^2
=2\!\int\rho\,\partial_t\rho
=\!-\int\rho(\mathbf{u}\cdot\nabla\rho)\,dx
+\!\int\rho(\partial_t\rho+\mu\Delta\rho)\,dx.
\end{equation}
For the $H^1$ norm,
\begin{equation}\label{eqn:pH1}
\begin{aligned}
\partial_t\|\nabla\rho\|_{L^2}^2
&=2\!\int \nabla\rho\cdot\nabla(\partial_t\rho)
=-2\!\int\partial_t\rho\,\Delta\rho \\
&=\frac{1}{2\mu}\!\int(\mathbf{u}\cdot\nabla\rho)^2\,dx
-
\frac{1}{2\mu}\!\int(\partial_t\rho+\mu\Delta\rho)^2\,dx.
\end{aligned}
\end{equation}

Substituting \eqref{eqn:H1}, \eqref{eqn:L2}, and \eqref{eqn:pH1} into \eqref{eqn:HoL} yields
\begin{equation}\label{eqn:comb}
\begin{aligned}
\partial_tH(t)
&=
\frac{1}{2\mu\|\rho\|_{L^2}^2}
\Bigl(
\int(\mathbf{u}\cdot\nabla\rho)^2\,dx
-
\int(\partial_t\rho+\mu\Delta\rho)^2\,dx
\Bigr) \\
&\quad
-\frac{1}{2\mu\|\rho\|_{L^2}^4}
\Bigl[
\Big(\!\int\rho\,\mathbf{u}\cdot\nabla\rho\Big)^2
-
\Big(\!\int\rho(\partial_t\rho+\mu\Delta\rho)\Big)^2
\Bigr].
\end{aligned}
\end{equation}

By Cauchy--Schwarz,
\[
\Bigl|\!\int\rho(\partial_t\rho+\mu\Delta\rho)\,dx\Bigr|^2
\le\|\rho\|_{L^2}^2
\!\int(\partial_t\rho+\mu\Delta\rho)^2\,dx.
\]
Therefore the negative term in \eqref{eqn:comb} cancels, and we obtain
\[
\partial_tH(t)
\le
\frac{\|\mathbf{u}\|^2_{L^\infty}}{2\mu}\,H(t).
\]
Gronwall's inequality now gives $H(t)\le C_1e^{C_1t}$, implying \eqref{eqn:ratio}.
\end{proof}

\medskip

\noindent\textbf{Conclusion of the proof of Theorem \ref{thm:para}.}
Multiplying \eqref{eqn:para} by $\rho$ and integrating, using $\nabla\cdot\mathbf{u}=0$, we obtain
\begin{equation}\label{eqn:int}
\partial_t\|\rho\|_{L^2}^2
=-2\mu\|\nabla\rho\|_{L^2}^2.
\end{equation}

Using Lemma~\ref{lem:ratio}, we have
\[
\|\nabla\rho(t,\cdot)\|_{L^2}^2
\le
C_1 e^{C_1t}\|\rho(t,\cdot)\|_{L^2}^2.
\]
Substituting this into \eqref{eqn:int} gives
\[
\partial_t\|\rho\|_{L^2}^2
\ge
-2\mu C_1 e^{C_1t}\|\rho\|_{L^2}^2.
\]
Another application of Grönwall’s inequality yields
\[
\|\rho(t,\cdot)\|_{L^2}
\ge
|\rho(T_0,\cdot)|_{L^2}\exp\!\bigl(-2\mu e^{C_1(t-T_0)}\bigr), \text{for all $t\geq T_0$.}
\]
which is precisely \eqref{eqn:double}.
As $\|\rho(t,\cdot)\|_{L^2}$ decreases in time, by choosing $C_2=\|\rho(T_0,\cdot)\|_{L^2}\exp\!\bigl(-2\mu e^{-C_1T_0}\bigr)$, we have \eqref{eqn:double}.

\end{proof}

\bibliographystyle{plain}
	\bibliography{citation}
\end{document}